\documentclass{article}
\usepackage{amsmath,amssymb}
\newtheorem{df}{Definition}[section]
\newtheorem{thm}[df]{Theorem}
\newtheorem{lem}[df]{Lemma}
\title{Uncertainty product of the spherical Abel-Poisson wavelet}
\author{I. Iglewska-Nowak\footnote{West Pomeranian University of Technology, Szczecin, School of Mathematics, al. Pias\-t\'ow 17, 70--310 Szczecin, Poland}}

\begin{document}

\maketitle

\bibliographystyle{amsplain}

\begin{abstract}In the present paper, the variances in the space and momentum domains of the spherical Abel-Poisson wavelet, as well as the limit of the uncertainty product for $\rho\to0$, where $\rho$ is the scale parameter, are computed. The values of these quantities coincide with a certain accuracy with those of the spherical Poisson wavelet with $\tfrac{1}{2}$ substituted for the order parameter.
\end{abstract}

\begin{bfseries}Key Words and Phrases:\end{bfseries} uncertainty product, time--frequency localization, Abel--Poisson wavelet\\
\begin{bfseries}2010 MSC:\end{bfseries} 42C40

\section{Introduction}

Similarly as in physics (Heisenberg's uncertainty principle), there exist several uncertainty principles in mathematics. An uncertainty constant of a function is a measure for the trade-off between the spatial and frequency localization. For functions over the sphere, the \emph{uncertainty product} was introduced by Narcowich and Ward in~\cite{NW96}, see also~\cite{IIN16USW} for a short history of the notion and a discussion about various nomenclature.

The present research is devoted to spherical wavelets (for the theory of wavelet transforms see~\cite{IIN15WT}). The uncertainty product can be interpreted as one of indicators of 'goodness' of wavelets, because well-localized wavelets yield a well-localized wavelet transform, i.e., such that 'fuzziness' of the transform reflects mostly the 'fuzziness' of the analyzed signal. Some general results concerning the uncertainty product of a wide class of wavelets in limit $\rho\to0$, where $\rho$ denotes the scale parameter, can be found in~\cite{IIN16USW}. The case of Poisson wavelets~$g_\rho^m$, $m\in\mathbb N$ \cite{HI07,IIN15PW} -- a wavelet family which is the most popular one for applications because of their explicit representation as well as existence of discrete frames \cite{IH10,IIN16WF} -- is investigated in~\cite{IIN17UPW}, and a bound for the uncertainty constant of the spherical Gauss-Weierstrass wavelet is computed in~\cite{IIN17UGW}.

In the present paper we investigate another wavelet family, namely the Abel-Poisson wavelet. The so-called variances in the space and momentum domains are computed explicitly and exactly. The uncertainty product of a function is the square root of the variances. In the case of the Abel-Poisson wavelet, its limit for $\rho\to0$ is finite. Up to a certain order, the values of the variances and the uncertainty product coincide with those for Poisson wavelets \cite{IIN17UPW} with $\tfrac{1}{2}$ substituted for the order parameter~$m$ (Poisson wavelets are defined for $m\geq1$).

The paper is organized as follows. After an introduction of the necessary notions and statements in Section~\ref{sec:sphere}, the main result of the paper, Theorem~\ref{thm:uncertaintyAP} is proven in Section~\ref{sec:uncertaintyAP}.

\section{Preliminaries}\label{sec:sphere}

Let $\mathcal{S}^n$ denote the $n$--dimensional unit sphere in $(n+1)$--dimensional Euclidean space~$\mathbb{R}^{n+1}$ with spherical variables $(\vartheta,\vartheta_2,\dots,\vartheta_{n-1},\varphi)$. Integrable zonal (rota\-tion-invariant with respect to the $x_1$-axis) functions over the sphere have the Gegenbauer expansion
$$
f(x)=\sum_{l=0}^\infty\widehat f(l)\,C_l^\lambda(\cos\vartheta)
$$
with Gegenbauer coefficients
$$
\widehat f(l)=c(l,\lambda)\int_{-1}^1 f(t)\,C_l^\lambda(t)\left(1-t^2\right)^{\lambda-1/2}dt,
$$
where $\lambda$ is an index related to the space dimension by
$$
\lambda=\frac{n-1}{2}
$$
and $c$ is a constant that depends on~$l$ and~$\lambda$. $C_l^\lambda$, $l\in\mathbb N_0$, are the Gegenbauer polynomials of order~$\lambda\in\mathbb R$ and degree~$l\in\mathbb{N}_0$. Unless it leads to misunderstandings, zonal spherical functions will be identified with those of the first variable, i.e., we write $f(x)$ interchangeable with $f(\cos\vartheta)$.

The variances in the space and momentum domains of a $\mathcal C^2(\mathcal S^n)$--function~$f$ with $\int_{\mathcal S^n}x\,|f(x)|^2\,d\sigma(x)\ne0$ are given by \cite{nLF03}
$$
\text{var}_S(f)=\left(\frac{\int_{\mathcal S^n}|f(x)|^2\,d\sigma(x)}{\int_{\mathcal S^n}x\,|f(x)|^2\,d\sigma(x)}\right)^2-1
$$
and
$$
\text{var}_M(f)=-\frac{\int_{\mathcal S^n}\Delta^\ast f(x)\cdot \bar f(x)\,d\sigma(x)}{\int_{\mathcal S^n}|f(x)|^2\,d\sigma(x)},
$$
where $\Delta^\ast$ is the Laplace--Beltrami operator on~$\mathcal S^n$. The quantity
$$
U(f)=\sqrt{\text{var}_S(f)}\cdot\sqrt{\text{var}_M(f)}
$$
is called the uncertainty product of~$f$.

The uncertainty product of zonal functions may be computed from their Gegenbauer coefficients \cite[Lemma~4.2]{IIN16MR} and according to the spherical uncertainty principle it is bounded from below by $\frac{n}{2}$, see \cite{NW96,RV97} and \cite[formula (4.37)]{GG04a}, \cite[formula (12)]{GG04b}.

\begin{lem}\label{lem:varS_varM} Let a zonal  square integrable and continuously differentiable function over $\mathcal S^n$ be given by its Gegenbauer expansion
$$
f(\cos\vartheta)=\sum_{l=0}^\infty\widehat f(l)\,\mathcal C_l^\lambda(\cos\vartheta).
$$
Its variances in space and momentum domain are equal to
\begin{align}
\text{var}_S(f)&=\left(\frac{\sum_{l=0}^\infty\frac{\lambda}{l+\lambda}\,\binom{l+2\lambda-1}{l}\,|\widehat f(l)|^2}{\sum_{l=0}^\infty\binom{l+2\lambda}{l}\,
   \frac{\lambda^2\left[\overline{\widehat f(l)}\,\widehat f(l+1)+\widehat f(l)\,\overline{\widehat f(l+1)}\right]}{(l+\lambda)(l+\lambda+1)}}\right)^2-1,\label{eq:varS}\\
\text{var}_M(f)&=\frac{\sum_{l=1}^\infty\frac{ l\lambda(l+2\lambda)}{l+\lambda}\,\binom{l+2\lambda-1}{l}\,|\widehat f(l)|^2}
   {\sum_{l=0}^\infty\frac{\lambda}{l+\lambda}\,\binom{l+2\lambda-1}{l}\,|\widehat f(l)|^2},\label{eq:varM}
\end{align}
whenever the series are convergent.
\end{lem}

\begin{thm}For $f\in\mathcal L^2(\mathcal S^n)\cap\mathcal C^1(\mathcal S^n)$, $U(f)\geq\frac{n}{2}$.
\end{thm}

The Abel-Poisson wavelet (with respect to the measure~$\alpha(\rho)=\frac{1}{\rho}$) is given by
$$
\Psi_\rho^A(x)=\sum_{l=0}^\infty\frac{\lambda+l}{\lambda}\sqrt{2l\rho}\,e^{-l\rho}\,\mathcal C_l^\lambda(\cos\vartheta),
$$
where $\rho\in\mathbb R_+$ denotes the scale parameter, see \cite[Section~3]{IIN15WT}.

\section{The uncertainty product of the Abel-Poisson wavelet}\label{sec:uncertaintyAP}

\begin{thm}\label{thm:uncertaintyAP}For the Abel-Poisson wavelet, the variances in space and momentum domain are given by
\begin{equation}\label{eq:varS_AP}\begin{split}
\text{var}_S(\Psi_\rho^A)=&\left\{\rho^n+2\left[2(n-1)(1-e^{-2\rho})^n\,\alpha(\rho)\,\rho^2-(4n^2-6n+3)\,\rho^n\right]\,e^{2\rho}\right.\\
&\left.-\left[4(n-1)(1-e^{-2\rho})^n\,\alpha(\rho)\,\rho^2+(4n-5)\,\rho^n\right]\,e^{4\rho}\right\}^{-2}\\
&\cdot\left\{16\,(n-1)^2\,\left[n-1+(n+1)\,e^{2\rho}\right]^2\,\rho^{2n}\,e^{2\rho}\right\}-1,
\end{split}\end{equation}
where~$\alpha$ is a bounded function, and
\begin{equation}\label{eq:varM_AP}
\text{var}_M(\Psi_\rho^A)=\frac{n\,(n+1)\left[n+(n+3)\,e^{2\rho}+e^{4\rho}\right]e^{2\rho}}{\left[n-1+(n+1)\,e^{2\rho}\right](e^{2\rho}-1)^2}.
\end{equation}
The uncertainty product equals in limit $\rho\to0$
\begin{equation}\label{eq:limitUP}
\lim_{\rho\to0}U(\Psi_\rho^A)=\frac{1}{2}\,\sqrt\frac{(n+1)(n+2)(n^2-3n+3)}{n(n-1)}.
\end{equation}
\end{thm}

\begin{bfseries}Proof. \end{bfseries}Substituting
$$
\widehat{\Psi_\rho^A}(l)=\frac{\lambda+l}{\lambda}\sqrt{2\rho l}\,e^{-\rho l}
$$
to the expressions~\eqref{eq:varS} and~\eqref{eq:varM} we obtain
\begin{align*}
\text{var}_S(\Psi_\rho^A)&=\left(\frac{\sum_{l=1}^\infty\frac{l+\lambda}{\lambda}\,\binom{l+2\lambda-1}{l}\,le^{-2\rho l}}
   {\sum_{l=1}^\infty\binom{l+2\lambda}{l}\cdot2\sqrt{l(l+1)}\,e^{-\rho(2l+1)}}\right)^2-1,\\
\text{var}_M(\Psi_\rho^A)&=\frac{\sum_{l=1}^\infty\frac{l(l+\lambda)(l+2\lambda)}{\lambda}\,\binom{l+2\lambda-1}{l}le^{-2\rho l}}
   {\sum_{l=1}^\infty\frac{l+\lambda}{\lambda}\,\binom{l+2\lambda-1}{l}\,le^{-2\rho l}}.
\end{align*}
Set
$$
S_{n,m}(\rho)=\sum_{l=1}^\infty\binom{l+2\lambda-1}{l}\,l^m\,e^{-2\rho l}.
$$
Then,
\begin{align}
\text{var}_S(\Psi_\rho^A)&=\left[\frac{e^\rho\cdot A(\rho)}{B(\rho)}\right]^2-1,\label{eq:varSPsi}\\
\text{var}_M(\Psi_\rho^A)&=\frac{C(\rho)}{A(\rho)}\label{eq:varMPsi}
\end{align}
for
\begin{align}
A(\rho)&=\frac{1}{\lambda}\,S_{n,2}(\rho)+S_{n,1}(\rho),\label{eq:Arho}\\
C(\rho)&=\frac{1}{\lambda}\,S_{n,4}(\rho)+3S_{n,3}(\rho)+2\lambda S_{n,2}(\rho),\label{eq:Crho}
\end{align}
and
$$
B(\rho)=\sum_{l=0}^\infty\frac{l+2\lambda}{\lambda}\binom{l+2\lambda-1}{l}\,\sqrt{l(l+1)}\,e^{-2\rho l}.
$$
In order to estimate $B(\rho)$ by $S_{n,m}(\rho)$ note that for $l\in\mathbb N$
$$
l+\frac{1}{2}-\frac{1}{8l}\leq\sqrt{l(l+1)}\leq l+\frac{1}{2}-\frac{1}{8l}+\frac{1}{16l^2}.
$$
Thus,
$$
\frac{l+2\lambda}{\lambda}\,\sqrt{l(l+1)}\geq\frac{l^2}{\lambda}+\left(\frac{1}{2\lambda}+2\right)l+\left(1-\frac{1}{8\lambda}\right)-\frac{1}{4l}
$$
and
\begin{align*}
\frac{l+2\lambda}{\lambda}\,\sqrt{l(l+1)}&\leq\frac{l}{\lambda}\left(l+\frac{1}{2}-\frac{1}{8l}+\frac{1}{16l^2}\right)+2\left(l+\frac{1}{2\lambda}\right)\\
&=\frac{l^2}{\lambda}+\left(\frac{1}{2\lambda}+2\right)l+\left(1-\frac{1}{8\lambda}\right)+\frac{1}{16\lambda l}.
\end{align*}
Consequently,
\begin{equation}\label{eq:Brho}
B(\rho)=\frac{1}{\lambda}\,S_{n,2}(\rho)+\left(\frac{1}{2\lambda}+2\right)S_{n,1}(\rho)+\left(1-\frac{1}{8\lambda}\right)S_{n,0}(\rho)+R(\rho)
\end{equation}
with
$$
|R(\rho)|\leq\frac{1}{4}\cdot S_{n,-1}(\rho).
$$
The values of $S_{n,1}$, $S_{n,2}$, $S_{n,3}$, and $S_{n,4}$ are obtained from
\begin{equation}\label{eq:Sn0}
S_{n,0}(\rho)=\frac{1}{(1-e^{-2\rho})^{n-1}}-1,
\end{equation}
compare the sixth equation on p.~352 in~\cite{IIN17UPW} (note that~$S_{n,0}$ in~\cite{IIN17UPW} is defined as a series starting from $l=0$), by the recurrence relation
\begin{equation}\label{eq:Snm_recurrence}
S_{n,m+1}(\rho)=-\frac{1}{2}S_{n,m}^\prime(\rho).
\end{equation}

In order to estimate~$S_{n,-1}(\rho)$ note that
$$
S_{n,-1}(\rho)=\sum_{l=1}^\infty\left[P_{n-3}(l)+\frac{(n-2)!}{l}\right]e^{-2\rho l},
$$
where~$P_{n-3}$ is a polynomial of degree~$n-3$. Thus, $S_{n,-1}$ can be expressed as a linear combination of $S_{2,m}$, $m=-1,0,\dots, n-3$,
\begin{equation}\label{eq:Sn-1_linear_combination}
S_{n,-1}(\rho)=\sum_{m=-1}^{n-3}\alpha_{n,m}\cdot S_{2,m}(\rho).
\end{equation}
It follows from~\eqref{eq:Sn0} and~\eqref{eq:Snm_recurrence} that for~$m\in\mathbb N_0$
\begin{equation}\label{eq:S2m}
S_{2,m}(\rho)=\frac{Q(n,e^{-2\rho})}{(1-e^{-2\rho})^{m+1}},
\end{equation}
where~$Q$ is a polynomial. Further,
$$
S_{2,-1}(\rho)=\sum_{l=1}^\infty\frac{e^{-2\rho l}}{l}.
$$
Since the function $t\mapsto\frac{e^{-2\rho t}}{t}$ is monotonously decreasing on $(0,\infty)$,
$$
\sum_{l=2}^\infty\frac{e^{-2\rho l}}{l}\leq\int_1^\infty\frac{e^{-2\rho t}}{t}\,dt.
$$
Thus,
\begin{equation}\label{eq:S2-1}
0\leq S_{2,-1}(\rho)\leq e^{-2\rho}+\int_1^\infty\frac{e^{-2\rho t}}{t}\,dt=e^{-2\rho}+\Gamma(0,2\rho).
\end{equation}
It follows from~\eqref{eq:Sn-1_linear_combination}, \eqref{eq:S2m}, \eqref{eq:S2-1}, and
$$
\frac{1}{(1-e^{-2\rho})^m}=\mathcal O\left(\frac{1}{\rho^m}\right)\qquad\text{for }\rho\to0
$$
that
$$
S_{n,-1}(\rho)=\mathcal O\left(\frac{1}{\rho^{n-2}}\right)\qquad\text{for }\rho\to0.
$$
Therefore, the rest term in~\eqref{eq:Brho} can be written as
$$
R(\rho)=\frac{\alpha(\rho)}{\rho^{n-2}}
$$
with a bounded function~$\alpha$. Formulae \eqref{eq:varS_AP} and \eqref{eq:varM_AP} are obtained from \eqref{eq:varSPsi} and \eqref{eq:varMPsi} by substituting \eqref{eq:Arho}, \eqref{eq:Brho}, and \eqref{eq:Crho} with the series~$S_{n,m}$, $m=1,2,3,4$, computed from \eqref{eq:Sn0} via \eqref{eq:Snm_recurrence}. In the computation, $S_{n,0}$ is replaced by~$(1-e^{-2\rho})^{1-n}$. The difference between the two expressions can be absorbed by the rest term~$R(\rho)$, and the replacement simplifies the calculations significantly.

\eqref{eq:limitUP} follows directly from \eqref{eq:varS_AP} and \eqref{eq:varM_AP}.\hfill$\Box$

The variances of the Abel-Poisson wavelet can be written as series
\begin{align*}
\text{var}_S(\Psi_\rho^A)&=\frac{n^2-3n+3}{n\,(n-1)}\cdot\rho^2-\frac{2n^2-4n+3+2^{n+1}\,n\,(n-1)\,\alpha(\rho)}{n^2\,(n-1)}\cdot\rho^3+\mathcal O(\rho^4),\\
\text{var}_M(\Psi_\rho^A)&=\frac{n^2+3n+2}{4\rho^2}+\frac{n^2-1}{2n\rho}+\mathcal O(1)
\end{align*}
for $\rho\to0$. Note that the first term of the variance in the space domain as well as the first two terms of the variance in the momentum domain are equal to those computed for Poisson wavelets \cite[Theorem~3.1]{IIN17UPW} with $\tfrac{1}{2}$ substituted for~$m$. Consequently, the limit of~$U(\Psi_\rho^A)$ for~$\rho\to0$ coincides with that for a Poisson wavelet with~$m=\tfrac{1}{2}$. Note that formally the Abel-Poisson wavelet is the Poisson wavelet of order~$\tfrac{1}{2}$. Therefore, the results concerning the variances and the uncertainty product are not surprising.

\end{document}